\documentclass[letterpaper, 10 pt, conference]{ieeeconf}
\IEEEoverridecommandlockouts
\usepackage{cite}
\usepackage{mathrsfs}
\usepackage{algorithmic}
\usepackage{graphicx}
\usepackage{textcomp}
\usepackage[dvipsnames]{xcolor}
\usepackage{svg}
\usepackage{float}
\usepackage{subcaption}
\usepackage[matmatrix, opt, reviewmark]{icpslab}
\usepackage{etoolbox}
\usepackage{booktabs}
\usepackage{multirow}
\usepackage{siunitx}
\usepackage{algorithm2e}

\newtheorem{remark}{Remark}

\newcommand{\MLPp}{\mathbf{HNN}_{\theta_p}}
\newcommand{\MLPd}{\mathbf{MLP}_{\theta_d}}
\newcommand{\numeq}{n_\mathrm{eq}}
\newcommand{\numineq}{n_\mathrm{in}}

\title{\LARGE \bf
HUANet: Hard-Constrained Unrolled ADMM for Constrained \\ Convex Optimization
}

\author{Trinh Tran$^\dagger$, Binh Nguyen$^\dagger$ and Truong X. Nghiem
\thanks{This material is based upon work supported by the U.S.\ National Science Foundation under Grants No.\;2449927 and No.\;2514584.}%
\thanks{All authors are with Department of Electrical and Computer Engineering,
University of Central Florida, Orlando, FL 32816, USA.
Emails: \{tr586225, thanhbinh.nguyen, truong.nghiem\}@ucf.edu
}
\thanks{$^\dagger$These authors contributed equally to this work.}
\thanks{Code is available at \url{https://github.com/nxt-lab/HUANet-code-public}.}
}

\begin{document}

    \maketitle
\thispagestyle{empty}
\pagestyle{empty}

\begin{abstract}
This paper presents HUANet, a constrained deep neural network architecture that unrolls the Alternating Direction Method of Multipliers (ADMM) into a trainable neural network for accelerating parametric constrained convex optimization.
Existing end-to-end learning methods operate as black-box mappings from parameters to solutions, often without explicitly incorporating optimality principles or guaranteeing constraint satisfaction. 
To address these limitations, HUANet embeds a hard-constrained neural network within each unrolled ADMM iteration, where a differentiable correction stage enforces the affine equalities of the primal subproblem.
Furthermore, we incorporate first-order optimality conditions into a self-supervised training loss to promote the convergence of the proposed unrolled algorithm. 
Extensive numerical experiments for benchmark optimization problems and a control application demonstrate and validate the effectiveness of HUANet in accelerating constrained convex optimization solving.

\end{abstract}

\section{Introduction}
Constrained convex optimization is an important class of mathematical programming problems that seeks to minimize a convex objective function subject to a set of convex constraints~\cite{bertsekas2014constrained}.
Such problems arise ubiquitously across a wide range of disciplines from engineering to finance, including robotics~\cite{nguyen2025spatially}, signal processing~\cite{gaines2018algorithms}, power grids~\cite{lin2012distributed}, transportation networks~\cite{mota2014distributed}, and portfolio optimization~\cite{jin2016constrained}.
These real-world optimization tasks are usually large-scale, involve varying problem parameters, and must satisfy safety-critical constraints.
Efficient algorithms for constrained convex optimization therefore remain an important research topic.
The key challenges for designing such an effective algorithm lie in the following performance criteria:
\textit{computational efficiency}, as the computational load increases substantially when problem dimensionality expands, and the solver must be re-executed for every new problem instance in parametric optimization problems;
\textit{feasibility}, as the solution must satisfy all constraints simultaneously;
and \textit{solution optimality}, as a solver should not sacrifice solution quality for speed.

Among many optimization algorithms available for solving such problems, the Alternating Direction Method of Multipliers (ADMM)~\cite{neal2011distributed, ADMM-Han2022} has emerged as a widely adopted first-order method due to its decomposability and ability to handle large-scale structured problems.
ADMM attempts to solve an optimization problem by decomposing it into smaller-sized subproblems, each of which is easier to solve.
While some ADMM subproblems admit closed-form analytical solutions, many practical problems do not.
Therefore, ADMM traditionally relies on numerical solvers, such as interior-point methods~\cite{nemirovski2008interior, goulart2024clarabel} or active-set methods~\cite{nocedal2006numerical}, to solve each subproblem.
However, these numerical solvers often require a large number of iterations for large-scale problems, creating a computational bottleneck that is particularly prohibitive for real-time applications.
Furthermore, ADMM suffers from slow tail convergence, demanding many iterations before reaching an optimal solution. 
These limitations motivate the development of a learning-based approach that can automatically adapt to varying problem instances, accelerate convergence, and guarantee constraint satisfaction.

Recently, learning-to-optimize (L2O) methods~\cite{chen2022learning} have emerged as a promising approach for solving parametric constrained optimization problems.
L2O models learn the solution mapping from problem parameters to solutions, which classical solvers do not exploit.
As a result, they can solve new problem instances significantly faster than classical methods.
Existing L2O methods can be broadly divided into two groups.
The first group consists of end-to-end surrogate models that replace the solver entirely by training a neural network to directly predict the solution~\cite{zamzam2020learning, wang2024fast}.
This approach uses supervised network training with training datasets of problem instances paired with optimal solutions from classical solvers.
While this approach has proven its power for unconstrained problems, it has struggled to perform well in constrained cases, as the predicted solutions are not guaranteed to satisfy the constraints at test time.
To address this shortcoming,~\cite{donti2021dc3} and~\cite{goertzen2026hardnet++} introduce predict-and-correct approaches that first use a neural network to predict a solution and then apply a correction step to project it onto the feasible set.
Such hybrid architectures can achieve higher feasibility rates than purely end-to-end approaches in some settings. 
However, as purely data-driven approaches, these deep learning methods lack interpretability since the networks learn a black-box mapping from inputs to outputs without any structural connection to the underlying optimization process. 

The second group uses machine learning to improve the performance of classical optimization solvers by embedding domain knowledge of the optimization algorithm into the network architecture.
One of the most prominent approaches in this group is algorithm unrolling~\cite{monga2021algorithm}.
Algorithm unrolling represents each iteration of a classical algorithm, such as the iterative shrinkage and thresholding algorithm~\cite{gregor2010learning}, gradient descent~\cite{andrychowicz2016learning}, and ADMM~\cite{yang2018admm, liu2025admm}, as one layer of the neural network.
These layers are then concatenated to form a deep neural network.
Compared to the first group, unrolling methods are more data efficient and more interpretable because they naturally inherit the prior structure and domain knowledge of the underlying algorithm. 
However, existing unrolling methods treat constraints as soft penalties, which provide no constraint feasibility guarantee for the solution.
In particular, existing unrolled ADMM methods~\cite{yang2018admm, liu2025admm}, which are the most relevant to our work, focus only on unconstrained problems, thus leaving a critical gap in applying the unrolled ADMM to general constrained optimization problems.
Additionally, despite inheriting the iterative 
structure of the underlying algorithm, unrolled networks do not guarantee convergence to an optimal solution. 
 
In this paper, we propose a Hard-Constrained Unrolled ADMM neural network (HUANet), an unrolled ADMM framework for accelerating the solution of parametric constrained convex optimization problems.
The main contributions of our approach are highlighted as follows:
\begin{itemize}
    \item \emph{Accelerating ADMM for constrained optimization with unrolled neural networks:} 
    We propose a hard-constrained neural network that enforces equality constraints via an end-layer correction to approximate the solution of the primal update with low computational cost at each ADMM iteration. 
    The proposed method provides a learnable and accelerated alternative to classical ADMM.
    
    \item \emph{Optimality-driven convergence for ADMM unrolling:}
     We develop a Karush–Kuhn–Tucker (KKT)-aware loss function that leverages terminal objective value, optimality residual terms and constraint penalties to guide network convergence.
     Alongside the hard-constrained primal network, we design a network that predicts the Lagrange multiplier of the primal subproblem. 
     The proposed loss function enables training in a self-supervised manner without requiring ground-truth solutions.

    \item \emph{Empirical validation on various constrained optimization problems}:
    We evaluate HUANet on constrained optimization benchmarks and a practical energy management problem, covering both smooth and nonsmooth constrained optimization settings.
    The results demonstrate that HUANet achieves significant computational speedups while maintaining high constraint satisfaction and solution accuracy compared with classical optimization solvers and existing learning-based methods.
\end{itemize}


\textbf{Notation:} Let $\bbR$, $\bbR_{\geq 0}$, and $\bbR_{>0}$ denote the sets of real numbers, nonnegative real numbers, and positive real numbers, respectively.
In Euclidean space $\bbR^n$, $\Vert x \Vert_2$ denotes the Euclidean norm of $x \in \bbR^n$.
For a vector $x \in \bbR^n$, $\calI_{\geq 0}(x)$ denotes the indicator function such that $\calI_{\geq 0}(x) = 0$ if $x \in \bbR_{\geq 0}^n$ and $\calI_{\geq 0}(x) = +\infty$  otherwise.
In addition, $\Pi_{\geq 0}(x)$ stands for the Euclidean projection of $x$ onto $\bbR_{\geq 0}^n$, that is, $\Pi_{\geq 0}(x) = \arg\min_{y\in \bbR_{\geq 0}^n}\Vert y - x \Vert_2^2$. 
Finally, $I_n$ denotes the identity matrix in $\bbR^{n\times n}$.

\section{Problem Formulation}\label{sec:problem-formulation}
Consider a parametric convex optimization problem of the following form 
\begin{subequations}
\label{eq:optimization_problem}
\begin{align}
    \minimize_x \; &f_\lambda(x),
    \\
    \text{s.t.} \; & A_\lambda x = b_\lambda, \label{eq:eq-constraint}
    \\
    & C_\lambda x \leq d_\lambda, \label{eq:ineq-constraint}
\end{align}
\end{subequations}
where 
$x \in \bbR^{n_x}$ is the decision variable,
$f_\lambda(x): \bbR^{n_x} \rightarrow \bbR \cup\{+\infty\}$ 
is a proper, closed, convex, and twice differentiable function on the interior of its domain,
and $A_\lambda \in \bbR^{\numeq \times n_x}, C_\lambda \in \bbR^{\numineq \times n_x}$, $b_\lambda \in \bbR^{\numeq}$, and $d_\lambda  \in \bbR^{\numineq}$.
Here, $\numeq$ and $\numineq$ denote the numbers of equality and inequality constraints, respectively.
In this paper, $f_\lambda$, $A_\lambda$,  $b_\lambda$, $C_\lambda$, and $d_\lambda$ are parameterized by the parameter vector $\lambda \in \Lambda \subset \bbR^{n_\lambda}$. 

By introducing a slack variable $s \in \bbR^{\numineq}$ and an auxiliary variable $w \in \bbR^{\numineq}$, the optimization problem~\eqref{eq:optimization_problem} can be rewritten equivalently as
\begin{subequations}
    \begin{align}
        \minimize_{x, s, w} \; &f_\lambda(x) + \mathcal{I}_{\geq 0}(w),
        \\
        \text{s.t.} \; &A_\lambda x=b_\lambda,
        \\
        &C_\lambda x + s = d_\lambda,
        \\
        & s = w.
    \end{align}
\end{subequations}
The problem can then be solved by the ADMM~\cite{neal2011distributed}, which comprises the following iterations:
\begin{subequations} \label{eq:admm}
\begin{align}
\label{eq:primal-update}
x^{k+1}, s^{k+1} &= 
\argmin_{x, s} \, f_\lambda(x) \!+\! \frac{\rho}{2}\norm{s \!-\! w^k \!+\! \rho^{-1}v^k}^2_2, 
\\
&\hspace{3em} \text{s.t.} \hspace{0.7em} A_\lambda x = b_\lambda, \nonumber
\\
& \hspace{5em}
C_\lambda x + s = d_\lambda,  \nonumber
\\
w^{k+1} &= \Pi_{\geq 0}(s^{k+1} + \rho^{-1}v^{k}), \label{eq:aux-update}
\\
v^{k+1} &= v^k + \rho(s^{k+1} - w^{k+1}) \label{eq:dual-update},
\end{align}
\end{subequations}
where $k = 0, 1, \ldots$ 
denotes the iteration index, $v \in \bbR^{\numineq}$ is the dual variable, and $\rho > 0$ is the step-size parameter.

However, in certain applications, such as model predictive control with a complicated cost function, a long horizon, 
or stringent real-time requirements, ADMM algorithm~\eqref{eq:admm} might be too computationally expensive.
Our \textbf{objective} is to develop a learning-based ADMM framework for accelerating the solution of the parametric optimization problem~\eqref{eq:optimization_problem}, while maintaining acceptable optimality gaps and constraint violations.

\begin{figure}[!tb]
    \centering
    \includegraphics[
    width=\linewidth,
    trim=1 1 1 1,
    clip
]{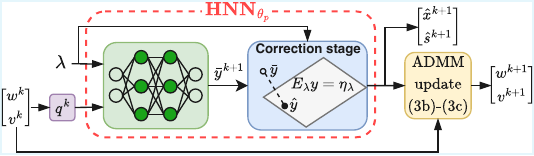}
    \caption{A neural layer $F_{\theta_p}$ of HUANet for an ADMM iteration at step $k$.     \label{fig:diag-L2O-a}}
\end{figure}

\section{Proposed HUANet Framework}\label{sec:proposed-huanet}

In this section, we present HUANet, a deep unrolled ADMM framework designed to efficiently solve the constrained optimization problem~\eqref{eq:optimization_problem} across varying problem parameters $\lambda$. HUANet unfolds the ADMM iterations in~\eqref{eq:admm} into $N$ sequential neural layers, where each layer corresponds to an ADMM iteration and executes an identical neural mapping $F_{\theta_p}$.
\subsection{HUANet Architecture}
\label{subsec:HUANet}
At each iteration $k$, the ADMM algorithm~\eqref{eq:admm} executes three sequential update steps, including the primal update~\eqref{eq:primal-update}, the auxiliary variable update~\eqref{eq:aux-update}, and the dual variable update~\eqref{eq:dual-update}.
Among these, the primal update is the most computationally expensive step, as it requires solving a constrained optimization subproblem.

Conventional end-to-end learning methods can approximate the optimal solution of the primal update.
They typically treat the inequality constraints as soft penalties, but do not guarantee constraint satisfaction in the primal update.

HUANet consists of $N$ neural layers, each defined by $F_{\theta_p}$, arranged sequentially and sharing the same learnable parameters $\theta_p$.
The architecture of each neural layer $F_{\theta_p}$ is illustrated in Fig.~\ref{fig:diag-L2O-a}, which includes a hard-constrained neural network, denoted by $\MLPp$.
The $\MLPp$ consists of a \textit{multi-layer perceptron (MLP)} that generates a prior estimate $(\bar x^{k+1}, \bar s^{k+1})$ of the primal–slack pair $(x^{k+1}, s^{k+1})$ without enforcing the equality constraints, followed by a \textit{correction stage}.
The outputs of the MLP are subsequently passed through this correction stage to enforce the equality constraints, as illustrated in Fig.~\ref{fig:diag-L2O-a}.
The correction stage acts as an equality feasibility corrector, maps $(\bar{x}^{k+1}, \bar{s}^{k+1})$ onto the affine subspace defined by the equality constraints, and produces a posterior estimate $(\hat x^{k+1}, \hat s^{k+1})$ of the primal–slack pair.

The \textit{MLP} is designed to learn the optimal solution of the primal update~\eqref{eq:primal-update}.
The inputs of the MLP are the vector $q^k  = w^k - \rho^{-1}v^k$ and the problem parameter $\lambda$.
The primal update can be rewritten as
\begin{subequations}
\begin{align}
\label{eq:primal-rewrite}
    x^{k+1}, s^{k+1} = \argmin_{x,s} \; &f_\lambda(x) 
+ \frac{\rho}{2}\left\|s - q^k\right\|_2^2, 
\\
\text{s.t.} \; & A_\lambda x = b_\lambda,\;C_\lambda x + s = d_\lambda.  \label{eq:primal-rewrite-cstr}
\end{align}
\end{subequations}
The MLP acts as a universal function approximator that learns the mapping $(q^k, \lambda) \mapsto \bar{y}^{k+1}=[(\bar{x}^{k+1})^\top,(\bar{s}^{k+1})^\top]^\top$, producing unconstrained estimates of the primal and slack variables prior to projection.
The MLP employs differentiable activation functions to enable efficient gradient-based optimization during training.

The \textit{correction stage} ensures that the output $\hat y^{k+1}$ of $\MLPp$ satisfies the equality constraints.
A common method is to project $\bar y^{k+1}$ onto the affine constraint set to obtain $\hat y^{k+1}$ by solving the following optimization problem
\begin{equation}
    \hat{y}^{k+1} = \argmin_{y} \frac{1}{2} \|y - \bar{y}^{k+1}\|_2^2,
    \quad \text{s.t.} \quad E_\lambda y = \eta_\lambda,
    \label{eq:projection-layer}
\end{equation}
where
$E_\lambda = \mm{
     A_\lambda, 0;
     C_\lambda, I_{\numineq}},
    \quad
    \eta_\lambda =\mm{
      b_\lambda;
      d_\lambda}\text.
$
If $E_\lambda$ has full row rank, 
the closed-form projection is
\begin{equation}
   \hat{y}^{k+1} = \bar{y}^{k+1} + E_\lambda^+ 
   (\eta_\lambda - E_\lambda \bar{y}^{k+1})\text,
   \label{eq:kkt-proj-closed-form}
\end{equation}
where $E_\lambda^+$ denotes the pseudo-inverse of $E_\lambda$.
Our method exploits the sparsity of $E_\lambda$ to compute $\hat{y}^{k+1}$ by solving the KKT conditions of \eqref{eq:projection-layer} with $LDL^\top$ factorization as implemented in OSQP and SCS~\cite{stellato2020osqp,scs2015}.

After obtaining $\hat{y}^{k+1}$, the ADMM variables $w^{k+1}$ and $v^{k+1}$ are updated according to~\eqref{eq:aux-update} and~\eqref{eq:dual-update}, respectively.

The neural network $\MLPp(q^k, \lambda)$, together with the ADMM updates~\eqref{eq:aux-update} and~\eqref{eq:dual-update} (by substituting $\hat s^{k+1}$ for $s^{k+1}$), constitutes a neural layer $F_{\theta_p}$, as illustrated in Fig.~\ref{fig:diag-L2O-a}.
It can be seen that the projected primal-slack pair $(\hat{x}^{k+1}, \hat{s}^{k+1})$ satisfies all equality constraints exactly, \ie both constraints $A_\lambda \hat{x}^{k+1} = b_\lambda$ and $C_\lambda \hat{x}^{k+1} + \hat{s}^{k+1} = d_\lambda$ always hold.
The HUANet framework is then constructed by stacking the neural layers $F_{\theta_p}$ sequentially.

From the diagram of $F_{\theta_p}$ in Fig.~\ref{fig:diag-L2O-a}, HUANet can be interpreted as a discrete-time, time-invariant dynamical system
\begin{subequations}
\label{eq:sys}%
\begin{align}
\mathbf{x}^{k+1} &= F_{\theta_p}(\mathbf{x}^k, \lambda),\quad\mathbf{x}^0 = \mm{w^0; v^0}, \label{eq:sys-1}
\\
\hat y^{k+1} &= \MLPp\left(\mm{I_{\numineq},\; -\rho^{-1} I_{\numineq}}\mathbf{x}^k, \lambda \right)\text,
\end{align}
\end{subequations}
where 
$\mathbf{x}^{k} = \mm{w^k; v^k}$, $\lambda$, and $\hat y^{k+1}$  are the state, the input, and the output of~\eqref{eq:sys}, respectively.
From~\eqref{eq:sys}, we can unroll the dynamics~\eqref{eq:sys-1} to accelerate the computation of HUANet.

It is important to note that in Fig.~\ref{fig:diag-L2O-a}, $N$ is the number of unrolled layers used only during training, whereas $K$ denotes the number of iterations performed at test time.

\begin{figure}[!tb]
    \centering
    \includegraphics[
    width=\linewidth,
    trim=1 1 1 1,
    clip
]{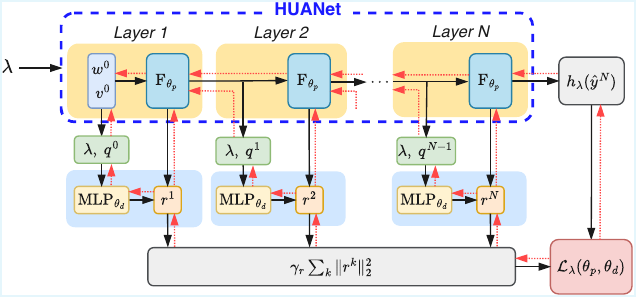}
    \caption{Training process of HUANet. Black and red arrows represent the forward and backward passes, respectively.}
    \label{fig:diag-L2O-b}
\end{figure}

\subsection{Self-supervised training} \label{subsec:admm-huanet}
This section describes the training procedure of HUANet for generating a high-accuracy feasible solution $\hat{x}^N$ for a given problem instance $\lambda$.
The training workflow, illustrated in Fig.~\ref{fig:diag-L2O-b}, consists of a forward pass and a backward pass.
Before detailing these passes, we first establish the optimality conditions that guide the training of HUANet.

\textit{KKT optimality conditions}: $\MLPp$ already guarantees the primal feasibility of~\eqref{eq:primal-rewrite}; however, the primal feasibility alone is insufficient to satisfy the optimality conditions.
According to~\cite[Chapter 5]{boyd2004convex}, sufficient optimality conditions for the primal update require both primal feasibility and stationarity conditions.
In this paper, the optimality conditions of the primal update are incorporated into the training of HUANet.
The stationarity conditions for subproblem~\eqref{eq:primal-rewrite} are
\begin{align}
    &\nabla f_\lambda(\hat{x}^{k+1}) + A_\lambda^\top z^{k+1} + C_\lambda^\top \beta^{k+1} = 0
    \label{eq:kkt-x},
    \\
    &\rho(\hat{s}^{k+1} - q^k) + \beta^{k+1} = 0,
    \label{eq:kkt-s}
\end{align}
where $z$ and $\beta$ denote the Lagrange multipliers associated with the constraints~\eqref{eq:primal-rewrite-cstr}.
Substituting~\eqref{eq:kkt-s} into~\eqref{eq:kkt-x} leads to
\begin{equation}
\label{eq:kkt-residual}
    r^{k+1} = \nabla f_\lambda(\hat{x}^{k+1}) + A_\lambda^\top z^{k+1} \!+
    \! \rho C_\lambda^\top(q^k \!-\! \hat{s}^{k+1}) = 0, 
\end{equation}
where $r^{k+1}$ is the KKT residual at iteration $k+1$.
The resulting residual $r^{k+1}$ serves as an augmentation term for the soft constraints in the training loss.
When $r^{k+1} = 0$, the  stationarity condition is satisfied, implying that the $(k+1)$-th layer of HUANet can mimic an ADMM iteration. 
Since $(\hat{x}^{k+1}, \hat{s}^{k+1})$ is itself determined by $(q^k, \lambda)$ via $\MLPp$, we treat $z^{k+1}$ as an implicit function of $(q^k, \lambda)$,  which captures its dependence on $(q^k, \lambda)$.
Accordingly, we use an additional MLP, denoted by $\MLPd$ and shown in Fig.~\ref{fig:diag-L2O-b}, to approximate the dual variable $z^{k+1} \approx \MLPd(q^k, \lambda)$.
The $\MLPd$ plays a crucial role in training HUANet, as it provides a learned dual variable $z^{k+1}$ for the primal update, which, together with the learned primal variable $\hat x^{k+1}$, enforces the optimality conditions. 
The KKT-augmented enforcement promotes convergence of HUANet, a property that is often overlooked in existing methods for unrolling conventional optimization algorithms.

We train HUANet in a self-supervised manner by minimizing a loss function that penalizes both the original optimization objective $f$ and violations of the KKT optimality condition~\eqref{eq:kkt-residual}.
The network is unrolled for $N$ ADMM iterations with all parameters trained jointly using the AdamW optimizer~\cite{loshchilov2017decoupled}.
For notational convenience, we index the unrolled iterations by $k = 1, \ldots, N$
hereafter.
For a batch of $S$ problem instances $\{\lambda^{(i)}\}_{i=1}^{S}$ sampled from a given distribution $\mathcal{D}$, the training loss is defined as
\begin{equation}
     \calL_\lambda({\theta_p}, \theta_d) = \frac{1}{S} \sum_{i=1}^{S}  \left(
     h_\lambda (\hat{y}_i^N) + 
     \frac{\gamma_r}{2} \sum_{k=1}^N \Vert r_i^k \Vert_2^2 \right),
    \label{eq:loss-function}
\end{equation}
where $h_\lambda(\hat{y}_i^N) = f_\lambda\bigl(\hat{x}_i^N\bigr) + \gamma_s \Vert \max(0, -\hat s_i^N)\Vert_2^2$ stands for the terminal cost,
which combines objective optimality with an inequality-violation penalty at the last neural layer,
${\theta_p}$ and $\theta_d$ denote trainable parameters of $\MLPp$ and $\MLPd$, $\gamma_s, \gamma_r > 0$ are weighting hyperparameters,
$x^k_i, s^k_i,$ and $r^k_i$ stand for the primal variable, slack variable, and KKT residual at layer $k$ of the $i$-th instance $\lambda^{(i)}$, respectively.
In the loss function~\eqref{eq:loss-function}, the first term enforces optimality and constraint satisfaction at the final layer, while the second term promotes satisfaction of the KKT conditions to facilitate convergence to the optimal solution.

\begin{remark}
The KKT-augmented loss function~\eqref{eq:loss-function} offers two practical advantages over existing unrolled methods.
First, since we consider first-order optimality conditions at each unrolled step, there is no need to accumulate the objective across all layers as in~\cite{andrychowicz2016learning}.
This can improve the accuracy and convergence of HUANet.

Second,~\eqref{eq:loss-function} defines a self-supervised loss function that does not require labeled data (computed optimal solutions) for the problem instances $\{\lambda^{(i)}\}_{i=1}^S$.
\end{remark}

\textit{Forward Pass:}
The complete forward pass of HUANet proceeds as follows.
Starting from the initial states $(w^0, v^0)$, with the input $\lambda^{(i)}$, each layer carries out the following sequential steps. 
The intermediate vector $q_i^k$ is computed from the current ADMM variables and, together with the problem parameter $\lambda^{(i)}$, is fed into $\MLPp$, which generates constrained estimates of the primal–slack pair.
The auxiliary and dual variables $w^{k+1}_i$ and $v^{k+1}_i$ are then updated via~\eqref{eq:aux-update} and~\eqref{eq:dual-update}.
At each layer, the KKT residual $r^{k+1}_i$ is computed via~\eqref{eq:kkt-residual}, and accumulated into the training loss. 
After $N$ layers, the total training loss $\calL_{\lambda}({\theta_p}, \theta_d)$ is computed over a batch of $S$ instances via~\eqref{eq:loss-function}.
The learned dual variables $\hat z^{k+1}_i = \MLPd(q_i^k, \lambda_i)$ are predicted by an MLP with trainable parameters $\theta_d$ and act as relaxation variables in the KKT residual of the training loss.

\textit{Backward Pass:} 
We update learnable parameters $({\theta_p}, \theta_d)$ of HUANet using first-order optimization methods, in which the gradient of $\calL_\lambda(\theta_p, \theta_d)$ can be computed by back-propagating through the unrolled layers. 
The derivatives of $\calL_\lambda$ with respect to the neural network parameters ${\theta_p}$ and ${\theta_d}$ are given by
\begin{align}
\pderiv{\calL_\lambda}{{\theta_p}} &= \frac{1}{S}\sum_{i=1}^{S} \left( \frac{dh_{\lambda}}{{d\theta_p}} + \gamma_r \sum_{k=1}^N (r_i^k)^\top \frac{dr_i^k}{d\theta_p}\right),
\label{eq:grad-L-theta}
\\
\pderiv{\calL_\lambda}{\theta_d} &= \frac{\gamma_r}{S} \sum_{i=1}^{S}  \sum_{k=1}^N (r_i^k)^\top \pderiv{r_i^k}{\theta_d}.
\label{eq:grad-L-alpha} 
\end{align}
By using the chain rule to compute the gradient terms of~\eqref{eq:grad-L-theta} and~\eqref{eq:grad-L-alpha}, we obtain
\begin{align}
\frac{d h_\lambda}{d {\theta_p}}
&=
\frac{\partial h_\lambda}{\partial \hat{y}_i^N}
\frac{\partial\hat{y}_i^N}{\partial \bar{y}_i^N}
\left(
\frac{\partial \bar{y}_i^N}{\partial\theta_p}
+
\frac{\partial \bar{y}_i^N}{\partial q_i^{N-1}}
\frac{d q_i^{N-1}}{d {\theta_p}}
\right),
\label{eq:df_dth}
\\
\frac{d r_i^k}{d{\theta_p}}
&=
\frac{\partial r_i^k}{\partial \hat{y}_i^k}
\frac{\partial\hat{y}_i^k}{\partial \bar{y}_i^k}
\left(
\frac{\partial \bar{y}_i^k}{\partial\theta_p}
+
\frac{\partial \bar{y}_i^k}{\partial q_i^{k-1}}
\frac{d q_i^{k-1}}{d{\theta_p}}
\right)
\nonumber\\
&\quad+
\left(
\frac{\partial r_i^k}{\partial q_i^{k-1}}
+
\frac{\partial r_i^k}{\partial z_i^k}
\frac{\partial z_i^k}{\partial q_i^{k-1}}
\right)
\frac{d q_i^{k-1}}{d{\theta_p}},
\label{eq:dr_dthp}
\\
\frac{\partial r_i^k}{\partial{\theta_d}}
&=
\frac{\partial r_i^k}{\partial z_i^k}
\frac{\partial z_i^k}{\partial{\theta_d}}.
\label{eq:dr_dthd}
\end{align}
In~\eqref{eq:df_dth}, \eqref{eq:dr_dthp}, and~\eqref{eq:dr_dthd},
$\frac{\partial \bar{y}_i^k}{\partial q_i^{k-1}}$,
$\frac{\partial \bar{y}_i^k}{\partial\theta_p}$,
$\frac{\partial z_i^k}{\partial q_i^{k-1}}$, and
$\frac{\partial z_i^k}{\partial\theta_d}$
can be computed via back-propagation through the $\MLPp$ and $\MLPd$,
respectively.
Furthermore,
\[
\frac{\partial r_i^k}{\partial \hat{x}_i^k}
= \nabla^2 f_\lambda(\hat x_i^k),\qquad
\frac{\partial r_i^k}{\partial \hat{s}_i^k}
= -\rho C_\lambda^\top,
\]
\[
\frac{\partial r_i^k}{\partial q_i^{k-1}}
= \rho C_\lambda^\top,\qquad
\frac{\partial r_i^k}{\partial z_i^k}
= A_\lambda^\top,
\]
as obtained from~\eqref{eq:kkt-residual}, and
from~\eqref{eq:kkt-proj-closed-form},
\begin{align}
\frac{\partial\hat{y}_i^k}{\partial \bar{y}_i^k}
=
I_{n_x+n_\mathrm{in}}-E_\lambda^+E_\lambda.
\end{align}

It should be noted that
$\frac{d q_i^{k-1}}{d{\theta_p}}$
represents the recurrent dependence induced by the preceding unrolled
ADMM iterations, while the remaining Jacobian terms can be computed
from variables within a single layer of HUANet.

\begin{remark} \label{rm:ext}
During training, the quantities $f(x_i^k)$, $\nabla f(x_i^k)$, and $\nabla^2 f_\lambda(x_i^k)$ must be well-defined for all $x_i^k$, which may not hold if $x_i^k$ lies outside the domain of $f$.
If the correction stage employs a projection operator to enforce the equality constraints in~\eqref{eq:primal-rewrite}, these quantities may not be well-defined at the projected points, which prevents the training process from being completed.
To address this issue, we can replace the correction stage with a \textit{feasibility stage} implemented as an optimization layer at the output of $\MLPp$.
Specifically, let $X$ denote the set of all points $x$ for which $f(x)$, $\nabla f(x)$, and $\nabla^2 f_\lambda(x)$ are well-defined. We first determine a feasible point $\hat {x}^k$ satisfying $\,\hat {x}^k \in X,\; A_\lambda \hat {x}^k = b_\lambda$. Then, we define $\hat {s}^k = d_\lambda - C_\lambda \hat {x}^k$. Accordingly, the output of the correction stage is given by $\hat {y}^k = \mm{(\hat {x}^k)^\top, (\hat {s}^k)^\top}^\top$.
\end{remark}

\begin{remark}
When there is no inequality constraint ($\numineq = 0$) in~\eqref{eq:optimization_problem}, HUANet has only one neural layer ($N=1$), which becomes the hard-constrained neural network $\MLPp$. In this case, $\MLPp$ and $\MLPd$ each take only $\lambda$ as input, and there are no ADMM updates.
\end{remark}

\section{Experimental Results} \label{sec:results}

To evaluate the performance of our proposed HUANet framework, we conduct numerical experiments for the quadratic programming (QP) problem, the entropy maximization problem, and the energy management problem~\cite{nguyen2026leaf}.

\subsection{Experimental Setup}

\textbf{Baselines}:
We compare HUANet with solvers like OSQP~\cite{stellato2020osqp}, Clarabel~\cite{goulart2024clarabel}, SCS~\cite{scs2015}, DC3~\cite{donti2021dc3}, and vanilla ADMM with Clarabel used to solve its subproblems.
Solver termination thresholds are adjusted to ensure comparable solution accuracies.

\textbf{Experimental setups}: All neural networks consist of two hidden layers, each with 512 neurons, across all optimization problems.
For each problem, a set of $20{,}000$ samples of $\lambda_i$ is generated, with $16{,}000$ used for training, $2{,}000$ for validation, and $2{,}000$ for testing.
For all problems, we set the numbers of training iterations $N$ and testing iterations $K$ within $[20,\;40]$, and allow a maximum of $500$ iterations for the ADMM baseline.

To ensure a fair comparison of computational time, HUANet inference times and baseline solver execution times are measured in wall-clock time on an Intel i9-285K CPU with 128 GB RAM.
HUANet is trained on an NVIDIA RTX 5090 GPU with 32 GB memory.

\textbf{Metrics}: We evaluate methods using optimality gap, equality violation, inequality violation, and solve time.
Constraint violations are reported as average (mean) and worst-case (max) infeasibility across test samples.

\subsection{Quadratic programming}

\subsubsection{Problem Setup}

This example validates the performance and accuracy of HUANet via 
parametric QP problems in the following form
\begin{subequations}
\begin{align*}
    \min_x \; &\frac{1}{2}x^\top Q x + p_\lambda^\top x,
    \\
    \text{s.t.} \; & Ax = b_\lambda,
    \\
    & Cx \leq d_\lambda,
\end{align*}
\end{subequations}
where $\lambda \in \bbR^{n_x}$, 
$(\lambda)_i \sim \calU(-1, 1)$ is the optimization parameter, 
$\calU(-1, 1)$ is the uniform distribution over $[-1, 1]$, 
$Q = F_q^\top F_q + I_{n_x}$ is a positive definite matrix,
$F_q \in \mathbb{R}^{n_x \times n_x}$ with
$(F_q)_{ij} \sim \mathcal{N}(0,1)$,
$p_\lambda = {\bf 1}_{n_x} + \lambda$,
$A\in \bbR^{\numeq\times n_x}$ with $(A)_{ij} \sim \calN(0,1)$,
and $C\in \bbR^{\numineq\times n_x}$ with $(C)_{ij} \;\sim \calN(0,1)$. 
In addition, 
$b_\lambda = - A Q^{-1} p_\lambda + \epsilon_b \in \bbR^{\numeq} $ and $d_\lambda = - C Q^{-1} p_\lambda + \epsilon_d \in \bbR^{\numineq}$ with 
$(\epsilon_b)_i \sim \calU(0, 0.1)$ and 
$(\epsilon_d)_i \sim \calU(0, 0.1)$.

\subsubsection{Results}

\begin{figure}[!tb]
    \centering
    \includegraphics[width=\linewidth]{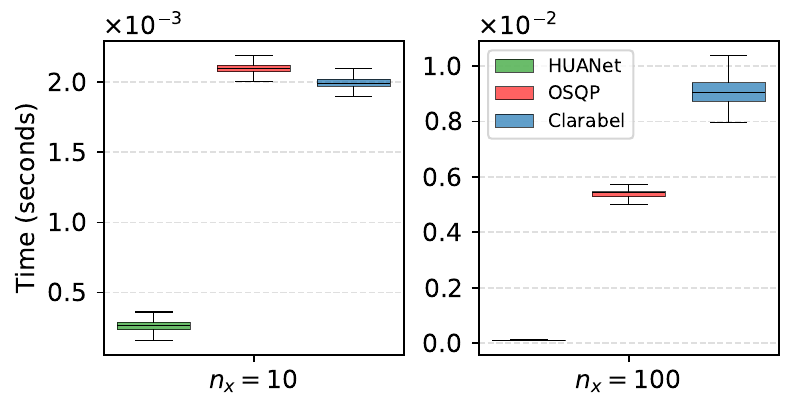}
    \caption{Run time comparison for the QP problem.}     \label{fig:qp_solve_time}
\end{figure}

\begin{figure}[!tb]
    \centering
    \begin{subfigure}[t]{0.48\linewidth}
        \centering
        \includegraphics[height=\linewidth]{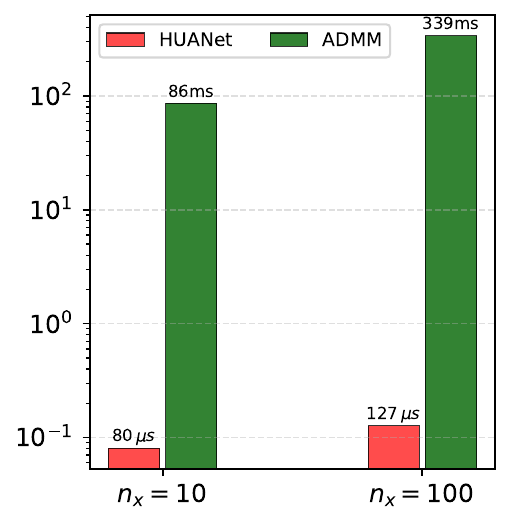}
        \caption{QP}
        \label{fig:qp_admm_time}
    \end{subfigure}
    \hfill
    \begin{subfigure}[t]{0.48\linewidth}
        \centering
        \includegraphics[height=\linewidth]{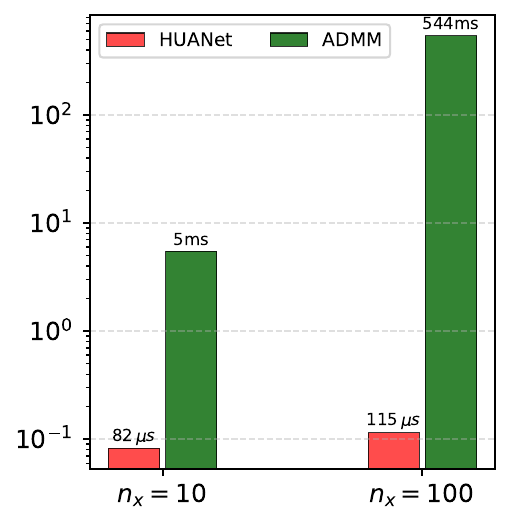}
        \caption{Entropy}
        \label{fig:entropy_admm_time}
    \end{subfigure}
    \caption{Average solve time of HUANet and ADMM for the QP and entropy problems.}
\end{figure}

Table~\ref{tab:results} shows that HUANet achieves a mean optimality gap below $0.06\%$ in both QP settings.
The correction stage enforces the affine equalities of the primal subproblem, resulting in equality violations on the order of $10^{-14}$ in numerical implementation.
Although inequality feasibility is not enforced by the correction stage, the terminal inequality penalty yields a mean inequality violation below $\num{3.6e-5}$ across both settings.
Fig.~\ref{fig:qp_solve_time} compares the execution time of HUANet against baselines.
For the low-dimensional problem ($n_x=10$), HUANet is approximately $7\times$ faster than solvers. 
As the problem dimensionality scales ($n_x=100$), the speedup becomes significantly more pronounced, with HUANet achieving an $88\times$ speedup over Clarabel and a $50\times$ speedup over OSQP.
Here, we formulate an intricate parametric QP problem where the problem parameters $\lambda$ appear in both the objective function and all constraints.
While this configuration increases the computational overhead for classical solvers, which must rebuild the optimization model for each new problem instance, HUANet avoids this bottleneck by leveraging the inference capabilities of a deep neural network.
As illustrated in Fig.~\ref{fig:qp_admm_time}, HUANet accelerates ADMM by unrolling the algorithm through the hard-constrained network.
Our approach achieves orders-of-magnitude speedups over vanilla ADMM, \eg a $1000\times$ speedup in the low-dimensional case and about a $2600\times$ speedup in the high-dimensional case.

\subsection{Entropy maximization}

\subsubsection{Problem Setup}

We consider the entropy maximization problem in~\cite[example 5.13]{boyd2004convex}
\begin{subequations}
\label{eq:maxEntropy}
\begin{align}
    \maximize_{x_i > 0}\; &-\sum_{i=1}^{n_x} x_i \log x_i 
    \\
    \text{s.t.\;} & {\bf 1}^\top x = 1, \label{eq:sum1}
    \\
    &Ax \leq b_\lambda, \label{eq:Axb}
\end{align}
\end{subequations}
where $x = [x_1, x_2,\dots, x_{n_x}]^\top \in \mathbb{R}^{n_x}$ is the decision vector representing a probability distribution, the matrix $A \in \mathbb{R}^{n_\mathrm{in} \times n_x}$ is sampled as $(A)_{ij} \sim \mathcal{N}(0,1)$, the problem parameter $\lambda \in \mathbb{R}^{n_\mathrm{in}}$ is given by $\lambda = A\upsilon + \epsilon$ with $\upsilon_i \sim \mathcal{U}(0,1)$ and $\epsilon_i \sim \mathcal{N}(0, 0.1)$, and we set $b_\lambda =\lambda$.
Here, we treat $b_\lambda$ as parameters of the optimization problem.

As mentioned in Remark~\ref{rm:ext}, to ensure that $f_i = x_i \log x_i$,
$f_i^\prime = \log x_i + 1$, and $f_i^{\prime\prime} = \frac{1}{x_i} $ are well-defined, the constraints $x_i > 0$ must be enforced in the correction stage. 
In this case, the correction stage first maps its output to a point $\hat{x}^{k+1}\in\bbR_{>0}^{n_x}$, scales this point to satisfy~\eqref{eq:sum1}, and then sets $\hat{s}^{k+1}=b_\lambda-A\hat{x}^{k+1}$.

\subsubsection{Results}

\begin{figure}[!tb]
    \centering
    \includegraphics[width=\linewidth]{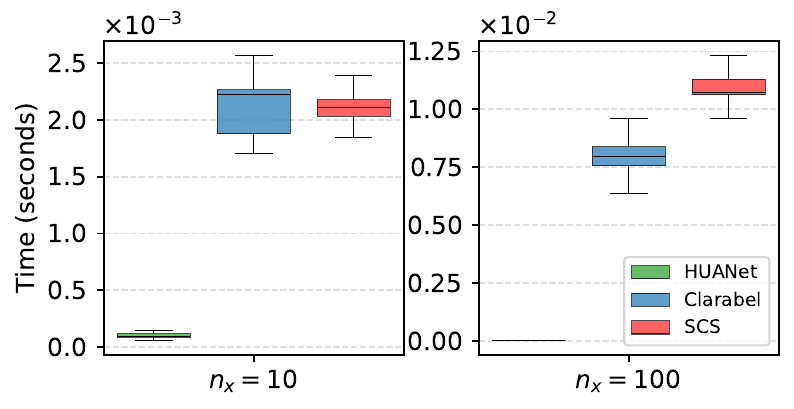}
    \caption{Run time comparison for the entropy problem.}
    \label{fig:entropy_solve_time}
\end{figure}

\begin{table*}[!t]
\centering
\caption{Performance across problem classes. QP and entropy results are
for HUANet; energy management results compare HUANet with DC3.
Entries report mean (maximum) over $1000$ test instances.}
\label{tab:results}
\begin{tabular}{lccrrrr}
\toprule
\textbf{Problem} & $n_x$ & $(\numeq,\numineq)$
& \multicolumn{1}{c}{\textbf{Optimality gap (\%)}}
& \multicolumn{1}{c}{\textbf{Eq. violation}}
& \multicolumn{1}{c}{\textbf{Ineq. violation}}
& \multicolumn{1}{c}{\textbf{Run time (s)}} \\
\midrule
QP
& 10 & (5,\,5)
& \num{4.56e-3}\;(\num{2.42e-2})
& \num{2.19e-15}\;(\num{3.55e-15})
& \num{6.99e-11}\;(\num{1.33e-8})
& \num{8.04e-5}\;(\num{2.28e-4}) \\
& 100 & (50,\,50)
& \num{5.94e-2}\;(\num{2.62e-1})
& \num{1.31e-14}\;(\num{2.49e-14})
& \num{3.55e-5}\;(\num{4.19e-4})
& \num{1.27e-4}\;(\num{2.97e-4}) \\
\midrule
Entropy
& 10 & (1,\,5)
& \num{2.83e-4}\;(\num{1.60e-2})
& \num{7.39e-17}\;(\num{4.44e-16})
& \num{2.25e-8}\;(\num{3.09e-6})
& \num{7.34e-5}\;(\num{1.33e-4}) \\
& 100 & (1,\,50)
& \num{6.18e-3}\;(\num{4.87e-2})
& \num{9.66e-17}\;(\num{4.45e-16})
& \num{9.56e-8}\;(\num{2.14e-6})
& \num{1.15e-4}\;(\num{1.96e-4}) \\
\midrule
Energy: HUANet
& 384 & (192, 481) 
& \num{5.64e-2}\;(\num{2.74e-1})
& \num{5.86e-15}\;(\num{7.49e-15})
& \num{3.81e-5}\;(\num{1.66e-3})
& \num{1.16e-3}\;(\num{1.33e-3}) \\
Energy: DC3
& 384 & (192, 481) 
& \num{4.05}\;(\num{5.67})
& \num{7.79e-15}\;(\num{1.42e-14})
& \num{1.54e-15}\;(\num{1.42e-14})
& \num{1.16e-4}\;(\num{3.40e-4}) \\
\bottomrule
\end{tabular}
\end{table*}
The results of HUANet on the entropy maximization problem are reported in Table~\ref{tab:results}.
The optimality gap remains below \num{6.18e-3}\% on average across both low- and high-dimensional settings.
Equality constraints are satisfied to machine precision, with a maximum violation of \num{4.45e-16}, while the maximum inequality violation is \num{3.09e-6}, comparable to the precision of the SCS solver.
The solve time comparison in Fig.~\ref{fig:entropy_solve_time} further demonstrates the computational advantage of HUANet.
For $n_x = 10$, HUANet is about $21\times$ faster than SCS and Clarabel.
These gaps widen significantly as problem size increases.
For $n_x = 100$, HUANet achieves a $235\times$ speedup over Clarabel and a $325\times$ speedup over SCS.
When the objective is nonquadratic, vanilla ADMM faces significant difficulties in solving the subproblems efficiently, whereas HUANet leverages the learned solution mapping to rapidly produce near-optimal feasible solutions.

\subsection{Energy management}

\subsubsection{Problem Setup}

Microgrid energy management requires repeatedly solving a constrained scheduling problem as demand, renewable-generation forecasts, and the measured battery state evolve.
In microgrid model predictive control, it is important for the controller to promptly update the energy schedule in response to changes in demand, renewable generation, and battery state.
We apply HUANet to the grid-connected microgrid in~\cite{nguyen2026leaf} to accelerate this recurring optimization while maintaining high solution accuracy.

The microgrid comprises photovoltaic (PV) generation, uncontrollable demand, a controllable load, and a battery energy storage system (BESS), which the controller schedules over a 24-hour horizon divided into $H=96$ intervals of duration $\Delta t=0.25\,\mathrm{h}$.
Let $\mathcal K=\{0,\ldots,H-1\}$ denote the interval indices.
At scheduling time $t$, the measured battery state of charge (SOC) is $\mathrm{SOC}_t$, and forecasts of uncontrollable demand $p_k^u$ and PV generation $p_k^{\mathrm{PV}}$ are available for every $k\in\mathcal K$.
Following~\cite{nguyen2026leaf}, we assume perfect forecasts over the horizon and define the net demand as $d_k=p_k^u-p_k^{\mathrm{PV}}$.
For each $k\in\mathcal K$, the decision variables are the battery power $p_k^B$, controllable-load consumption $p_k^c$, and grid exchange $m_k$.
The optimization also determines the SOC trajectory $\mathrm{SOC}_0,\ldots,\mathrm{SOC}_H$.
Positive $p_k^B$ denotes battery discharge, whereas negative $p_k^B$ denotes charging.
Similarly, $m_k>0$ denotes grid import, whereas $m_k<0$ denotes grid export.
The power balance $p_k^c+d_k=m_k+p_k^B$ equates total consumption with power supplied by the grid and battery.

The stage cost $J_k$ at interval $k$,
\begin{equation*}
J_k=r_{\mathrm{ec}}\Delta t\bigl(m_k+\alpha|p^B_k|\bigr)
+r_{\mathrm{op}}\max(m_k,0)+r_{\mathrm{df}}D(p_k^c),
\end{equation*}
includes the energy cost, the positive grid import charge, and the discomfort cost, respectively.
The coefficient $r_{\mathrm{ec}}>0$ is the energy price ($\$/\unit{\kWh}$), $r_{\mathrm{op}}\geq0$ weights positive grid import ($\$/\unit{\kW}$), and $r_{\mathrm{df}}>0$ is the discomfort-cost coefficient ($\$$).
The term $\alpha|p_k^B|$ in the energy cost accounts for battery cycling losses,
where $\alpha=(1-\eta)/(2\sqrt{\eta})$ is the loss coefficient and $\eta\in(0,1]$ is the battery round-trip efficiency.
Following~\cite{nguyen2026leaf}, the positive-import charge is accumulated at each stage.
The dimensionless discomfort function $D(\cdot)$ is defined as
\begin{equation*}
D(p)=
\begin{cases}
+\infty, & p\le 0, \\
a/p-1, & 0<p<a, \\
0, & p\ge a
\end{cases}
\end{equation*}
for a preferred controllable-load consumption $a>0$ ($\unit{kW}$).
For $p>0$, this function reduces to $D(p)=\max(a/p-1,0)$.
Then, the control objective is to minimize the total electricity and discomfort cost $J = \sum_{k=0}^{H-1} J_k$ over the prediction horizon while satisfying the power balance, battery dynamics, operating limits, and terminal-energy requirement.

The 
terms $|p^B_k|$, $\max(m_k,0)$, and $\max(a/p-1,0)$
admit exact epigraph reformulations, which introduce $3H$ auxiliary variables and $5H$ affine inequalities.
Because the HUANet training procedure requires a twice continuously differentiable objective function, we instead approximate them using the smooth convex functions
\begin{align*}
\psi_\delta(z) &= \sqrt{z^2+\delta^2}-\delta \approx |z|\text,\\
\phi_\kappa(z) &= \frac{\log(1+\exp(\kappa z))}{\kappa} \approx \max(z, 0)\text,
\end{align*}
with smoothing parameters $\delta>0$ and $\kappa>0$.
The stage cost is then approximated by
\begin{multline*}
    \tilde{J}_k = r_{\mathrm{ec}}\Delta t \left(m_k+\alpha\psi_\delta(p^B_k)\right)
    + r_{\mathrm{op}}\phi_{\kappa_m}(m_k) \\ + r_{\mathrm{df}}\phi_{\kappa_d}\left(\frac{a}{p_k^c}-1\right)\text.
\end{multline*}
We set $\delta=0.01\,\unit{kW}$, $\kappa_m=100$, 
and $\kappa_d = 10$. 
We also impose a fixed lower bound $p^c_{\min}>0$ to keep the controllable load away from the singularity of $a/p_k^c$ at zero.

The resulting scheduling problem is given by
\begin{align}
\minimize_{\mathrm{SOC},p^B,p^c,m}\;&
\tilde{J} = \sum_{k=0}^{H-1} \tilde{J}_k \label{eq:em-smooth}\\ 
\text{s.t.}\;&
\mathrm{SOC}_{k+1}
=\mathrm{SOC}_k-\frac{\Delta t}{b_{\mathrm{ess}}}p^B_k,\; k\in \mathcal{K},\notag \\ 
&p_k^c+d_k=m_k+p^B_k,\; k\in \mathcal{K}, \notag \\ 
&p^B_{\min}\le p^B_k\le p^B_{\max},\; k\in \mathcal{K}, \notag \\ 
&p^c_{\min}\leq p_k^c,\; k\in \mathcal{K},\notag \\ 
&\mathrm{SOC}_{\min}\le\mathrm{SOC}_k\le\mathrm{SOC}_{\max},\; k=0,\ldots,H,\notag \\ 
&\mathrm{SOC}_0=\mathrm{SOC}_t,
\quad \mathrm{SOC}_H\geq\overline{\mathrm{SOC}}, \notag\label{eq:SOC-limit-2}
\end{align}
where $b_{\mathrm{ess}}>0$ denotes the BESS energy capacity ($\unit{kWh}$),
$\mathrm{SOC}_t$ is the measured SOC at scheduling time $t$, and $\overline{\mathrm{SOC}}$ denotes the minimum terminal SOC.
To map \eqref{eq:em-smooth} to the standard form \eqref{eq:optimization_problem}, we stack the trajectories as $m=[m_0,\ldots,m_{H-1}]^\top$, 
$p^B=[p^B_0,\ldots,p^B_{H-1}]^\top$,
$p^c=[p^c_0,\ldots,p^c_{H-1}]^\top$,
and $\mathrm{SOC}=[\mathrm{SOC}_1,\ldots,\mathrm{SOC}_H]^\top$,
each in $\mathbb R^H$.
We then define
\begin{align*}
x&=[m^\top,(p^B)^\top,(p^c)^\top,\mathrm{SOC}^\top]^\top
\in\mathbb R^{4H},\\
\lambda&=\mathrm{SOC}_t.
\end{align*}
Denote $c_B=\Delta t/b_{\mathrm{ess}}$ and the matrix $L=I_H-
\bigl[\,0_{(H-1)\times1}\ \ I_{H-1}\,\bigr]^\top
\bigl[\,I_{H-1}\ \ 0_{(H-1)\times1}\,\bigr]
\in\mathbb R^{H\times H}$.
Let $e_1,e_H\in\mathbb R^{H}$ be the first and last canonical basis vectors, respectively.
Then 
the constraints \eqref{eq:eq-constraint} and \eqref{eq:ineq-constraint} have
\begin{equation*}
A_\lambda=
\begin{bmatrix}
0 & c_B I_H & 0 & L\\
I_H & I_H & -I_H & 0
\end{bmatrix},
\qquad
b_\lambda=
\begin{bmatrix}
\lambda e_1\\
d
\end{bmatrix},
\end{equation*}
\begingroup
\setlength{\arraycolsep}{2.5pt}
\begin{equation*}
C_\lambda=
\begin{bmatrix}
0 & I_H & 0 & 0\\
0 & -I_H & 0 & 0\\
0 & 0 & -I_H & 0\\
0 & 0 & 0 & I_H\\
0 & 0 & 0 & -I_H\\
0 & 0 & 0 & -e_H^\top
\end{bmatrix},
\;
d_\lambda=
\begin{bmatrix}
p^B_{\max}\mathbf 1_H\\
-p^B_{\min}\mathbf 1_H\\
-p^c_{\min}\mathbf 1_H\\
\mathrm{SOC}_{\max}\mathbf 1_H\\
-\mathrm{SOC}_{\min}\mathbf 1_H\\
-\overbar{\mathrm{SOC}}
\end{bmatrix}.
\end{equation*}
\endgroup
Here, $0_H\in\mathbb R^H$ is the zero vector, $\mathbf1_j\in\mathbb R^j$ is the all-ones vector, $I_j\in\mathbb R^{j\times j}$ is the identity matrix, and unlabeled zero blocks have the dimensions implied by the block rows and columns.
The physical parameters and cost weights follow~\cite{nguyen2026leaf}.
\subsubsection{Results}
\begin{figure}[!tb]
    \centering
    \includegraphics[width=\linewidth]{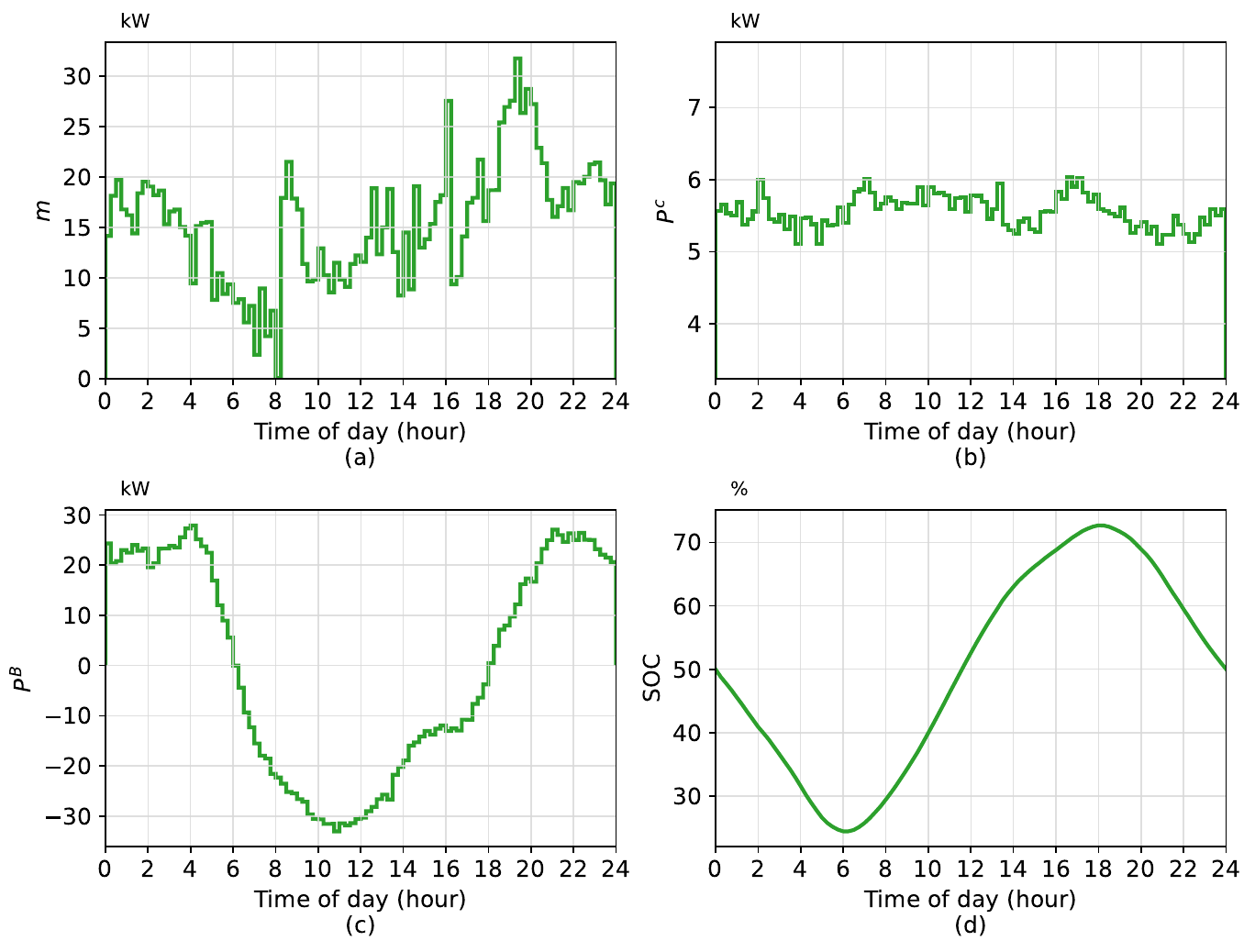}
    \caption{Day-ahead microgrid schedule predicted by HUANet over a 24-hour horizon.}     \label{fig:energy-performance}
\end{figure}

As shown in Table~\ref{tab:em_runtime}, HUANet is approximately $12\times$ faster than Clarabel in this experiment.
The difference relative to vanilla ADMM is larger: HUANet requires approximately $1\,\mathrm{ms}$ per instance, whereas vanilla ADMM requires $3\,\mathrm{s}$.
Such short latency is useful in receding-horizon operation, where the energy-management problem is solved repeatedly as demand and renewable-generation forecasts are updated.
Fig.~\ref{fig:energy-performance} illustrates the day-ahead schedule obtained by HUANet.
The controllable load remains within a narrow operating range, while the battery alternates between charging and discharging as net demand varies throughout the day.
The equality constraints are satisfied to machine precision, and the mean inequality violation is $\num{3.81e-5}$, as reported in Table~\ref{tab:results}.

Table~\ref{tab:results} compares the performance of HUANet and DC3.
While DC3 has a lower mean inference time ($0.116\,\mathrm{ms}$) than HUANet ($1.18\,\mathrm{ms}$), HUANet achieves drastically higher solution accuracy.
Specifically, it achieve an optimality gap of just $0.0564\%$ compared to $4.05\%$ for DC3, representing nearly two orders of magnitude improvement.
Thus, the comparison shows a tradeoff between DC3's solving time and HUANet's higher solution accuracy.
This comparison does not by itself establish safety or closed-loop performance.


\begin{table}[!tb]
    \centering
    \caption{Runtime comparison for the energy problem.}
    \label{tab:em_runtime}
    \footnotesize
    \setlength{\tabcolsep}{12pt}
    \begin{tabular}{lcc}
        \hline
        \textbf{Method} & \textbf{Run time (ms)} & \textbf{Speedup*} \\
        \hline
        HUANet       & $1.16 \pm 0.05$          & $12.58\times$ \\
        Clarabel     & $14.60 \pm 0.47$         & $1.00\times$  \\
        Vanilla ADMM & $3026.93 \pm 2669.44$    & $0.005\times$ \\
        \hline
    \end{tabular}
    \vspace{3pt}
    
    {\scriptsize
    *Speedup is measured relative to Clarabel.}
\end{table}

All the above experiments clearly showcase the effectiveness, scalability, and generalizability of HUANet for parametric constrained convex optimization.

\section{Conclusion}
\label{sec:conclusion}
In this paper, we introduced HUANet, a deep unrolling architecture for accelerating parametric constrained convex optimization by representing each ADMM iteration as a neural layer.
The proposed architecture enforces satisfaction of the equality constraints at every unrolled iteration via a correction stage in each neural layer.
By integrating the first-order optimality conditions directly into the training process, our framework facilitates convergence of the unrolled iterates toward optimal solutions.
Our architecture achieves high performance on various constrained convex problems while offering significant computational efficiency compared to vanilla ADMM and traditional solvers, particularly in high-dimensional scenarios.
For future work, we will embed both equality and inequality constraints into the network architecture, thus eliminating the reliance on explicit correction stages.
In addition, we will consider broader classes of convex problems beyond affine constraints, and ultimately extend the framework to nonconvex problems.

\bibliographystyle{IEEEtran}
\bibliography{references}

\end{document}